\documentclass[12pt]{article}

\usepackage[affil-it]{authblk}
      \usepackage{latexsym}
         \usepackage[reqno, namelimits, sumlimits]{amsmath}
         \usepackage{amssymb, amsfonts}
         \usepackage{amsthm}

 \newtheorem{theorem}{Theorem}[section]
 
 \newtheorem{definition}[theorem]{Definition}
 
 \newtheorem{lemma}[theorem]{Lemma}

 \newtheorem{cor}[theorem]{Corollary}
 \newtheorem{pro}[theorem]{Proposition}

\title{Local Regularity of Axisymmetric Solutions to the Navier-Stokes Equations}

\author{ G Seregin
  \thanks{seregin@maths.ox.ac.uk; }
  }
\affil{OxPDE, Mathematical Institute, University of Oxford, Oxford,UK, and 
St Petersburg Department of V A Steklov Mathematical Institute, Russia}
\date{ \today}

\begin{document}
\maketitle

\centerline{Dedicated to Vladimir  Gilelevich Mazya}

\begin{abstract} In the note, a local regularity condition for axisymmetric solutions to the non-stationary 3D Navier-Stokes equations is proven. It reads that  axially symmetric energy solutions to the Navier-Stokes equations have no Type I blowups.

\end{abstract}

\setcounter{equation}{0}
\section{Introduction }

The aim of the note is to discuss  potential singularities of axisymmetric solutions to the non-stationary 3D Navier-Stokes equations. Roughly speaking, we would like to show that if scale-invariant energy quantities of an axially symmetric solution are bounded then such a solution  is smooth. 
By definition, potential singularities with bounded scale-invariant energy quantities are called Type I blowups. It is important to notice that our result does not follow from the so-called $\varepsilon$-regularity theory, where regularity is coming from  smallness of those scale-invariant energy quantities.



Before stating and proving the main result of the note, see Theorem \ref{mainresult}, we are going to  remind basic notions from the mathematical theory of the Navier-Stokes equations.

For simplicity, let us consider the following Cauchy problem for the Navier-Stokes equations:
\begin{equation}\label{NSE}
\partial_tv+v\cdot\nabla v-\Delta v=-\nabla q,\qquad {\rm div}\,v=0	
\end{equation}
in $Q_+=\mathbb R^3\times ]0,\infty[$ and 
\begin{equation}\label{ini}
v|_{t=0}=u_0	
\end{equation}
in $\mathbb R^3$, where $u_0\in C^\infty_{0,0}(\mathbb R^3):=\{v\in C^\infty_{0}(\mathbb R^3):\,{\rm div}\,v=0\}$.
 
One of the main problems of the mathematical theory of viscous incompressible fluids is the global well-posedness of the Cauchy  problem (\ref{NSE}) and (\ref{ini}). A plausible approach (which, of course, is not unique) is to prove the global existence of  a solution  and then to prove its uniqueness.   It was done by J. Leray many years ago in \cite{JL1934}, who introduced the notion which is known now as  a weak Leray-Hopf solution. 

\begin{definition} \label{wLH}
A divergence free velocity field $v$ is a weak Leray-Hopf solution to the Cauchy problem (\ref{NSE}) and (\ref{ini}) if it has the following properties:

1. $v\in L_\infty(0,\infty;J)\cap L_2(0,\infty;J^1_2)$, where $J$ is the closure of $ C^\infty_{0,0}(\mathbb R^3)$ in $L_2(\mathbb R^3)$ and $J^1_2$ 	the closure of $ C^\infty_{0,0}(\mathbb R^3)$ with respect to the semi-norm $\Big(\int\limits_{\mathbb R^3}|\nabla w|^2dx\Big)^\frac 12$;

2. the function $t\mapsto\int\limits_{\mathbb R^3}v(x,t)\cdot w(x)dx$ is continuous on $[0,\infty[$ for any $w\in L_2(\mathbb R^3)$;

3. the Navier-Stokes equations is satisfied as the variational identity
$$\int\limits_{Q_+}\Big\{v\cdot\partial_tw+v\otimes v:\nabla w-\nabla v:\nabla w\Big\}dxdt=0$$ for any test vector-valued function $w\in C^\infty_0(Q_+)$ with ${\rm div}\,w=0$;

4. $\|v(\cdot,t)-u_0(\cdot)\|_{L_2(\mathbb R^3)}\to 0$ as $t\downarrow 0$;

5. global energy inequality 
$$\frac 12 \int\limits_{\mathbb R^3}|v(x,t)|^2dx+\int\limits^t_0\int\limits_{\mathbb R^3}|\nabla v|^2dxdt'\leq \frac 12 \int\limits_{\mathbb R^3}|u_0|^2dx$$ holds for all $t\geq 0$.
\end{definition}

 There is no information about the pressure in Definition \ref{wLH}.  But one can easily  recover the pressure by means of the linear theory and the following estimate for the pressure takes place 


$$
\|q\|_{L_\frac 32(Q_+)}\leq c\|v\|^2_{L_3(Q_+)}.$$

There is an important class of weak solutions, 
which is closely related to the uniqueness of solutions to the initial boundary value problems for the Navier-Stokes equations.
\begin{definition}\label{strong}
	Let $v$ be a weak Leray-Hopf solution to (\ref{NSE}) and (\ref{ini}). It is a strong solution to the Cauchy problem on the set $Q_T=\mathbb R^3\times ]0,T[$ if $\nabla v\in L_{2,\infty}(Q_T)$.
\end{definition}
The proposition, proved by J Leray in \cite{JL1934}, essentially reads the following. Assume that $v$ a weak Leray-Hopf solution to (\ref{NSE}) and (\ref{ini}). Then there exists a number $T\geq c\|\nabla u_0\|_{L_2(\mathbb R^3)}^{-4}$, where $c$ is an universal constant, such that $v$ is strong solution in $Q_T$.

Another important result proved by J Leray is the so-called weak-strong uniqueness. Assume that $v^1$ is another weak Leray-Hopf solution with the same initial data $u_0$, then $v^1=v$ on $Q_T$.

So, as it follows from the above statement, in order to prove uniqueness of weak Leray-Hopf solution on the interval $]0,T[$, it is enough to show that $\nabla v\in L_{2,\infty}(Q_T)$. In other words, the problem of unique solvability of the Cauchy problem in the energy class can be reduced to the problem of regularity of weak Leray-Hopf solutions.

As usual  in the theory of non-linear equations, we do not need to prove that $\nabla v\in L_{2,\infty}(Q_T)$ for some $T>0$. In fact, it is enough to prove a weaker regularity result and then the remaining part of the proof of regularity  follows from the linear theory. For example, one of the convenient spaces for such intermediate regularity is the space $L_\infty(Q_T)$.

So, the first time when a singularity occurs can be defined as follows:
$$\limsup\limits\limits_{t\uparrow T}\|v(\cdot,t)\|_{L_\infty(\mathbb R^3)}=\infty.$$

To study regularity of weak Leray-Hopf solutions by classical PDE's methods, we should mimic energy solutions on the local (in space-time) level. To this end,  the pressure should be involved into  considerations. The corresponding setting has been already discussed by Caffarelli-Kohn-Nirenberg, who have introduced the notion of suitable weak solutions to the Navier-Stokes equations, see \cite{CKN}, \cite{Sc}, and \cite{Sc1}.
\begin{definition}\label{sws}
Let $\omega\subset \mathbb R^3$ and $T_2>T_1$. $w$ and $r$ is a suitable weak solution to the Navier-Stokes in $Q_*=\omega\times ]T_1,T_2[$ if:

1. $w\in L_{2,\infty}(Q_*)$, $\nabla w\in L_2(Q_*)$, $r\in L_\frac 32(Q_*)$;

2. $w$ and $r$ satisfy the Navier-Stokes equations in the sense of distributions;

3. for a.a. $t\in [T_1,T_2]$, the local energy inequality
$$\int\limits_\omega\varphi(x,t)|w(x,t)|^2dx+2
\int\limits_{T_1}^t\int\limits_\omega\varphi|\nabla w|^2dxdt'\leq\int\limits_{T_1}^t\int\limits_\omega[|w|^2(\partial_t\varphi+\Delta \varphi)+$$
$$+w\cdot\nabla\varphi(|w|^2+2r)]dxdt'$$	holds for all non-negative $\varphi\in C^1_0(\omega\times ]T_1,T_2+(T_2-T_1)/2[).$
\end{definition}


In order to state a typical result of the  regularity theory for suitable weak solutions, introduce the notation for parabolic cylinders (balls): $Q(z_0,R)=B(x_0,R)\times ]t_0-R^2,t_0[$, $B(x_0,R)=\{x\in\mathbb R^3:\,\,|x-x_0|<R\}$, and $z_0=(x_0,t_0)$.
\begin{pro}\label{regsws}

1.	There are universal constants $\varepsilon$ and $c_0$ such that for any suitable weak solution $v$ and $q$ in $Q(z_0,R)$ satisfying the assumption
$$C(z_0,R)+D(z_0,R)<\varepsilon,$$
where
$$C(z_0,R)=\frac 1{R^2}\int\limits_{Q(z_0,R)}|v|^3dz,\qquad D(z_0,R)=\frac 1{R^2}\int\limits_{Q(z_0,R)}|q|^\frac 32dz,$$
the velocity field $v$ is H\"older continuos in $\overline Q(z_0,R/2)$ and 
$$\sup\limits_{zQ(z_0,R/2)}|v(z)|\leq \frac {c_0}R.$$

2. There is a universal constant $\varepsilon>0$ such that for any suitable weak solution $v$ and $q$ in $Q(z_0,R)$ satisfying the assumption
$$g(z_0):=\min\{\limsup\limits_{r\to0}E(z_0,r),\limsup\limits_{r\to0}A(z_0,r),\limsup\limits_{r\to0}C(z_0,r)\}<\varepsilon, $$
where
$$A(z_0,R)=\sup\limits_{t_0-R^2<t<t_0}\frac 1R\int\limits_{B(x_0,R)}|v(x,t)|^2dx,\qquad
E(z_0,R)=\frac 1R\int\limits_{Q(z_0,R)}|\nabla v|^2dz,$$
the point $z_0$ is a regular point of $v$, i.e., there exists $0<r\leq R$ such that $r\in L_\infty(Q(z_0,r))$.
\end{pro}
All the involved quantities are invariant with respect to the Navier-Stokes scaling 
$$v(x,t)\to \lambda v(\lambda x, \lambda^2t),\qquad q(x,t)\to \lambda^2 q(\lambda x, \lambda^2t).$$
We often call them energy scale invariant quantities.

Proposition \ref{regsws} describes the so-called  $\varepsilon$-regularity theory of suitable weak solutions. In this context,
it is interesting to verify whether or not a given weak Leray-Hopf solution has the property to be a suitable weak one in subdomains.
\begin{theorem}\label{turbulent}
	Among of all weak Leray-Hopf solutions to the Cauchy problem with the same initial data $u_0$, there exists at least one solution which is a suitable weak solution in any $Q(z_0,R)\subset Q_+$.\end{theorem}

	\begin{cor}\label{turbulent2} A weak Leray-Hopf solution, having the properties indicated in Theorem \ref{turbulent}, is, in fact, a turbulent solution, i.e., there exists a set  $S$ of full measure, containing zero, with the following property: for any $s\in S$, the inequality 
	$$\int\limits_{\mathbb R^3}|v(x,t)|^2dx+2\int\limits^t_s \int\limits_{\mathbb R^3}|\nabla v|^2dxdt'\leq \int\limits_{\mathbb R^3}|v(x,s)|^2dx	$$ holds for any $t\leq s$.		\end{cor}

Notion of turbulent solutions has been introduced in \cite{JL1934}.  


Finally, we would like to define different types of blowups (singularities).


Let $v$ and $q$ be a suitable weak solution to the Navier-Stokes equations in $Q(z_0,R)$. Obviously, if 
$$g(z_0)=\min\{\limsup\limits_{r\to0}E(z_0,r),\limsup\limits_{r\to0}A(z_0,r),\limsup\limits_{r\to0}C(z_0,r)\}>0,$$
then $z_0$ is a singular point.
\begin{definition}\label{deftypeI}
Let $z_0$ be a singular point. It is a Type I 	blowup if $g(z_0)<\infty$. The point $z_0$ is of  Type II if $g(z_0)=\infty$.\end{definition}
It is useful to notice that if $g(z_0)<\infty$, see \cite{Seregin2006}, then 
$$G(z_0)=\max\{\limsup\limits_{r\to0}E(z_0,r),\limsup\limits_{r\to0}A(z_0,r),\limsup\limits_{r\to0}C(z_0,r),$$$$\limsup\limits_{r\to0}D_0(z_0,r)\}<\infty.$$
Here, 
$$D_0(z_0,r)=\frac 1{r^2}\int\limits_{Q(z_0,r)}|q-[q]_{B(x_0,r)}|^\frac 32dz.$$

In fact, as it follows from \cite{SZ2007}, see Lemma 2.1 there, $D_0(z_0,r)$ can be replaced with 
$$D(z_0,r)=\frac 1{r^2}\int\limits_{Q(z_0,r)}|q|^\frac 32dz.$$

		\label{conmaincond}

\setcounter{equation}{0}
\section{Axial Symmetry and Blowups}
There are many papers on regularity of axially symmetric solutions. We cannot pretend to cite all good works in this direction. For example, 
let us mention papers: \cite{L1968}, \cite{UY1968},
\cite{LMNP1999}, \cite{NP2001}, \cite{CL2002}, \cite{LeiZhang2011}, \cite{ChenStrainYauTsai2009}, and \cite{ChenFangZhang2017}.

In this section, we assume that $\Omega$ is the unit cylinder centred at the origin, i.e.,
$$\Omega =\mathcal C=\{x=(x',x_3), \,x'=(x_1,x_2):\,|x'|<1,\,|x_3|<1\}.$$

In what follows, it will be convenient to replace balls $B(x_0,r)$ with cylinders $\mathcal C(x_0,r)=\{x=(x',x_3):\,|x'-x'_0|<r,\,|x_3-x_{03}|<r\}$. We let
$\mathcal C(r)=\mathcal C(0,r)$.

In our standing assumption, it is supposed that a suitable weak solution $v$ and $q$ to the Navier-Stokes equations in $Q=\mathcal C\times ]-1,0[$ is axially symmetric with respect to the axis $x_3$. The latter means the following: if we introduce the corresponding cylindrical coordinates $(\varrho,\varphi,x_3)$ and use the corresponding representation $v=v_\varrho e_\varrho+v_\varphi e_\varphi+v_3e_3$, then $v_{\varrho,\varphi}=v_{\varphi,\varphi}=v_{3,\varphi}=q_{,\varphi}=0$.

The main statement of the note reads that axisymmetric solutions have no  Type I singulaities.
 \begin{theorem}
\label{mainresult}
Assume that a pair $v$ and $q$	is axially symmetric suitable weak solution to the Navier-Stokes equations in $Q$ 
and the origin $z=0$ is a singular point of $v$. Then it is a Type II blowup.
\end{theorem}

\begin{proof} 
	Our fist observation  is related to
	 Lemma 3.3 in \cite{SS2009}. It is stated there that:
\begin{equation}
	\label{locMoser}
\sup\limits_{z\in Q(1/2)} \sigma(z)\leq C(M)\Big(\int\limits_{Q( 3/4)}|\sigma(z)|^\frac {10}3dz\Big)^\frac 3{10},
\end{equation}
where $\sigma :=\varrho v_\varphi=v_2x_1-v_1x_2$ and
$$M=\Big(\int\limits_{Q( 3/4)}|\overline v(z)|^\frac {10}3dz\Big)^\frac 3{10}, \qquad \overline v=v_re_r+v_3e_3.$$
However, estimate (\ref{locMoser}) has been proven in \cite{SS2009} under the additional assumption that there are no singular points if $t<0$. Let us try to modify the proof of Lemma 3.3 in \cite{SS2009} to a  general case of axially symmetric suitable weak solutions. Indeed, the basic part of the proof will remain to be Moser iterations.

So, let us show  that, under assumptions  of Theorem \ref{mainresult},  the left hand side of (\ref{locMoser}) is finite.  As in the proof of Lemma 3.3 of \cite{SS2009}, we hope to estimate $\omega^\frac {10}3 $ in a smaller domain by $\omega^\frac 52$ in larger domain, where $\omega=\sigma^{m}$ with a suitable number $m\geq 1$. To this end, let us try to understand how smooth  axially symmetric suitable weak solutions are. Denote by $S$ the set of singular points of $v$. It is well known that $S$ has 1D Hausdorff measure zero, $x'=0$ for any $z=(x,t)\in S$, and any spatial derivative of $v$ is H\"older continuous in $\mathcal C\times ]-1,0]\setminus S$. 

It is easy to check that $\sigma$ satisfies the 
equation 
\begin{equation}
	\label{eqsigma}
	\partial_t\sigma+\Big(v+2\frac{(x',0)}{|x'|^2}\Big)\cdot\nabla\sigma-\Delta\sigma=0
\end{equation}
in $Q\setminus(\{x'=0\}\times]-1,0[)$. Given $N>0$, set $\sigma_N(x,t)=\sigma (x,t)$ if $|\sigma(x,t)|\leq N$, $\sigma_N(x,t)= N$ if $\sigma>N$, and  $\sigma_N(x,t)=-N$ if $\sigma<-N$. Take two cut-off functions: the first of them is $\psi=\psi(x,t)$, vanishing in a neighbourhood of the parabolic boundary of $Q$ and the second one is $\phi=\phi(x')$, vanishing in a neighbourhood of the axis of symmetry,
 multiple the left hand side of equation (\ref{eqsigma}) by $\sigma_N^{2m-1}\psi^4\phi^2$ and integrate by parts.   As a result,  three different terms appear and they will be treated separately. For the first one, we have
$$\int\limits_{-1}^{t_*}\int\limits_\mathcal C\partial_t\sigma\sigma_N^{2m-1}\psi^4\phi^2dz=\frac 1{2m}\int\limits_{\mathcal C}\sigma_N^{2m}\psi^4\phi^2|_{(x,t_*)}dx+$$
$$+\int\limits_{\mathcal C}(\sigma -\sigma_N)\sigma_N^{2m-1}\psi^4\phi^2|_{(x,t_*)}dx-$$
$$-\frac 1{2m}\int\limits_{-1}^{t_*}\int\limits_{\mathcal C}\sigma_N^{2m}\partial_t\psi^4\phi^2dz-
\int\limits_{-1}^{t_*}\int\limits_{\mathcal C}(\sigma -\sigma_N)\sigma_N^{2m-1}\partial_t\psi^4\phi^2dz.
$$ 
Since the second term on the right hand side of the above identity is non-negative, it can be dropped out: 
$$\int\limits_{-1}^{t_*}\int\limits_\mathcal C\partial_t\sigma\sigma_N^{2m-1}\psi^4\phi^2dz\geq\frac 1{2m}\int\limits_{\mathcal C}\sigma_N^{2m}\psi^4\phi^2|_{(x,t_*)}dx-$$
$$-\frac 1{2m}\int\limits_{-1}^{t_*}\int\limits_{\mathcal C}\sigma_N^{2m}\partial_t\psi^4\phi^2dz-
\int\limits_{-1}^{t_*}\int\limits_{\mathcal C}(\sigma -\sigma_N)\sigma_N^{2m-1}\partial_t\psi^4\phi^2dz.$$
Denoting $b=2(x',0)|x'|^{-2}$,  transform
the second term as follows:
$$\int\limits_{-1}^{t_*}\int\limits_{\mathcal C}(v+b)\cdot\nabla \sigma\sigma_N^{2m-1}\psi^4\phi^2dz=
$$
$$=-\frac 1{2m}\int\limits_{-1}^{t_*}\int\limits_{\mathcal C}(v+b)\cdot\nabla (\psi^4\phi^2)\sigma_N^{2m}dz-$$$$-\int\limits_{-1}^{t_*}\int\limits_{\mathcal C}(v+b)\cdot\nabla (\psi^4\phi^2)(\sigma-\sigma_N)\sigma_N^{2m-1}dz
$$
Finally, for the third term, we have
$$-\int\limits_{-1}^{t_*}\int\limits_{\mathcal C}\Delta \sigma \sigma_N^{2m-1}\psi^4\phi^2dz=$$$$
=\frac {2m-1}{m^2}\int\limits_{-1}^{t_*}\int\limits_{\mathcal C}|\nabla \sigma_N^m|^2\psi^4\phi^2dz
-\frac 1{2m}\int\limits_{-1}^{t_*}\int\limits_{\mathcal C}\sigma_N^{2m}\Delta (\psi^4\phi^2)dz-$$$$-
\int\limits_{-1}^{t_*}\int\limits_{\mathcal C}(\sigma-\sigma_N)\cdot \sigma_N^{2m-1}\Delta (\psi^4\phi^2)dz. 
$$
Combining previous relationships, we find the following energy inequality
$$\frac 1{2m}\int\limits_{\mathcal C}\sigma_N^{2m}\psi^4\phi^2|_{(x,t_*)}dx+\frac {2m-1}{m^2}\int\limits_{-1}^{t_*}\int\limits_{\mathcal C}|\nabla \sigma_N^m|^2\psi^4\phi^2dz\leq$$$$
\leq\frac 1{2m}\int\limits_{-1}^{t_*}\int\limits_{\mathcal C}\sigma_N^{2m}(\partial_t(\psi^4\phi^2)+(\Delta \psi^4\phi^2))dz+$$$$+
\int\limits_{-1}^{t_*}\int\limits_{\mathcal C}(\sigma -\sigma_N)\sigma_N^{2m-1}(\partial_t(\psi^4\phi^2)+\Delta( \psi^4\phi^2))dz+$$
$$+\frac 1{2m}\int\limits_{-1}^{t_*}\int\limits_{\mathcal C}(v+b)\cdot\nabla (\psi^4\phi^2)\sigma_N^{2m}dz+$$$$+\int\limits_{-1}^{t_*}\int\limits_{\mathcal C}(v+b)\cdot\nabla (\psi^4\phi^2)(\sigma-\sigma_N)\sigma_N^{2m-1}dz.$$

Now,  selecting a special non-negative cut-off function $\phi$ so that $\psi(x')=0$ if $0<|x'|<\varepsilon/2$, $\psi(x')=1$ if $|x'|>\varepsilon$, and $|\nabla^k\phi|\leq c\varepsilon^{-k}$, $k=0,1,2$, let us see what happens if $\varepsilon\to 0$. We start with the two most  important terms:
$$I_1=\int\limits_{-1}^{t_*}\int\limits_{\mathcal C}|v|\psi^4\phi|\nabla\phi|(|\sigma|+|\sigma_N|)|\sigma_N|^{2m-1}dz
$$
and 
$$
I_2=\int\limits_{-1}^{t_*}\int\limits_{\mathcal C}|b|\psi^4\phi|\nabla\phi|(|\sigma|+|\sigma_N|)|\sigma_N|^{2m-1}dz.$$
As to $I_1$, it is easy to see
$$I_1\leq c\int\limits_{-1}^{t_*}2\pi\int\limits^1_{-1}\int\limits_{\varepsilon/2<\varrho<\varepsilon}|v|^2 N^{2m-1}\varrho d\varrho dx_3dt\to0
$$ as $\varepsilon\to0$. The second term is the most difficult one. Indeed,
$$I_2\leq c\int\limits_{-1}^{t_*}2\pi\int\limits^1_{-1}\int\limits_{\varepsilon/2<\varrho<\varepsilon}\frac 1\varepsilon N^{2m-1}|v_\varphi|\varrho d\varrho dx_3 dt\leq $$
$$\leq cN^{2m-1}\frac 1\varepsilon\Big(\int\limits_{-1}^{t_*}\int\limits^1_{-1}\int\limits_{\varepsilon/2<\varrho<\varepsilon}|v_\varphi|^2dz\Big)^\frac 12\Big(\int\limits_{-1}^{t_*}\int\limits^1_{-1}\int\limits_{\varepsilon/2<\varrho<\varepsilon}dz\Big)^\frac 12\to 0$$
as $\varepsilon\to0$.
 It remains to estimate the first two terms in the energy inequality. In the worst case scenario, we proceed as  follows:
$$I_3=cN^{2m-1}\frac 1{\varepsilon^2}\int\limits_{-1}^{t_*}\int\limits^1_{-1}\int\limits_{\varepsilon/2<\varrho<\varepsilon}|\sigma|dz=cN^{2m-1}\frac 1{\varepsilon}\int\limits_{-1}^{t_*}\int\limits^1_{-1}\int\limits_{\varepsilon/2<\varrho<\varepsilon}|v_\varphi|dz\to 0$$
as $\varepsilon\to0$, see bounds for $I_2$. So, passing to the limit in the energy inequality, we find 
$$\frac 1{2m}\int\limits_{\mathcal C}\sigma_N^{2m}\psi^4|_{(x,t_*)}dx+\frac {2m-1}{m^2}\int\limits_{-1}^{t_*}\int\limits_{\mathcal C}|\nabla \sigma_N^m|^2\psi^4dz\leq$$$$
\leq\frac 1{2m}\int\limits_{-1}^{t_*}\int\limits_{\mathcal C}\sigma_N^{2m}(\partial_t\psi^4+\Delta \psi^4)dz+$$$$+
\int\limits_{-1}^{t_*}\int\limits_{\mathcal C}(\sigma -\sigma_N)\sigma_N^{2m-1}(\partial_t\psi^4+\Delta \psi^4)dz+$$
$$+\frac 1{2m}\int\limits_{-1}^{t_*}\int\limits_{\mathcal C}(v+b)\cdot\nabla \psi^4\sigma_N^{2m}dz+$$$$+\int\limits_{-1}^{t_*}\int\limits_{\mathcal C}(v+b)\cdot\nabla \psi^4(\sigma-\sigma_N)\sigma_N^{2m-1}dz.$$

Now, pick up  a cut-off function exactly as in \cite{SS2009}, i.e., $\psi(x,t)=\Phi(x)\chi(t)$ so that
$0\leq \Phi\leq 1$, $\Phi=0$ outside $\mathcal C(r_1)$, $\Phi=1 $ in $\mathcal C(r)$,  $0\leq\chi\leq 1$, $\chi(t)=0$ if $t\leq -r_1^2$, $\chi(t)=1$ if  $t>-r^2$,  and 
$$\nabla^k \Phi\leq \frac c{|r_1-r|^k}, \quad k=1,2,\qquad |\partial_t\chi|\leq \frac c{(r_1-r)^2},
$$
where $1/2\leq r<r_1\leq 3/4$. Letting $\omega=\sigma^m$ and $\omega_N=\sigma^m_N$, we derive from the energy inequality the  basic estimate:
$$|\psi^2\omega_N|^2_{2,Q}:=\sup\limits_{-1<t<0}\|\psi^2\omega_N\|^2_{L_2(\mathcal C)}+\|\nabla(\psi^2\omega_N)\|^2_{L_2(Q)}\leq$$
\begin{equation}
	\label{energy1}
\leq \frac c{(r_1-r)^2}\int\limits_{Q(r_1)}\psi^2\omega\omega_Ndz+cJ_1+cJ_2,\end{equation}
where 
$$J_1=\frac 1{r_1-r}\int\limits_{Q(r_1)}|v||\omega ||\psi^2\omega_N|dz$$
and 
$$J_2=\frac 1{r_1-r}\int\limits_{Q(r_1)}\psi^3|b||\omega |^\frac 1m|\omega_N|^{2-\frac 1m}dz$$

Now, we split the proof into three  steps.

\textit{Step I} We assume that $m=4/3$, $r_1=3/4$, and $r=5/8$.
Then we have
$$\frac c{(r_1-r)^2}\int\limits_{Q(r_1)}\psi^2\omega\omega_Ndz\leq c\frac {r_1}{(r_1-r)^2}\Big(\int\limits_{Q(r_1)}|\omega|^\frac 52dz\Big)^\frac 45.$$
For the second term, we are going to exploit the H\"older inequality. Indeed,\begin{equation}\label{intermed}
	J_1\leq \frac 1{r_1-r}\Big(\int\limits_{Q_1(r_1)}|v|^\frac {10}3dz\Big)^\frac 3{10}\Big(\int\limits_{Q_1(r_1)}|\omega|^\frac {5}2dz\Big)^\frac 2{5}\Big(\int\limits_{Q_1(r_1)}|\psi^2\omega_N|^\frac {10}3dz\Big)^\frac 3{10}.	\end{equation}
Finally, taking into account that $\omega =(|x'||v_\varphi|)^\frac 43$, we estimate the third term:
$$J_2\leq c\int\limits_{Q_1(r_1)}\frac 1{|x'|}\omega^2dz\leq c\int\limits_{Q_1(r_1)}|v|^\frac 83dz\leq c
\Big(\int\limits_{Q(r_1)}|v|^\frac {10}3dz\Big)^\frac 45$$
Using the known multiplicative inequality  
\begin{equation}\label{mult}
\|\psi^2\omega_N\|_{L_{\frac {10}3}(Q(r_1))}\leq c|\psi^2\omega_N|_{2,Q(r_1)}	
\end{equation}
and then  Young inequality, we can pass to the limit as $N\to\infty$ and conclude that
\begin{equation}
	\label{step1}
	\int\limits_{Q(5/8)}|\sigma|^{(\frac {4}3)^2\frac 52}dz=\int\limits_{Q(5/8)}|\sigma|^{\frac {4}3\frac {10}3}dzd<\infty.
\end{equation}

\textit{Step II} Now,  assume that $m\geq (4/3)^2$, $r_1\leq 5/8$. The first two terms in energy inequality (\ref{energy1})
can be evaluated in the same way as in Step I. So, let us focus on the third term:
$$J_2\leq \frac 1{r_1-r}\int\limits_{Q(r_1)}\psi^3|b||\omega |^\frac 9{16}|\omega_N|^{\frac {23}{16}}dz\leq$$
$$\leq\frac 1{r_1-r}\int\limits_{Q(r_1)}|b||\omega |^\frac 9{16}|\psi^2\omega_N|^{\frac {23}{16}}dz $$
 For $s=40/21$,  it follows from the H\"older inequality that
$$J_2\leq \frac 1{r_1-r}\Big(\int\limits_{\mathcal C(r_1)}|b|^sdx\Big)^\frac 1s\int\limits^0_{-r_1^2} \Big(\int\limits_{\mathcal C(r_1)} |\omega|^\frac {s'9}{16} |\psi^2\omega_N|^\frac {s'23}{16}dx\Big)^\frac 1{s'},$$
where $s'=s/(s-1)=40/19$ and 
$$J_3=\Big(\int\limits_{\mathcal C(r_1)}|\omega|^\frac {s'9}{16}|\psi^2\omega_N|^\frac {s'23}{16}dx\Big)^\frac 1{s'}\leq $$$$\leq \Big(\int\limits_{\mathcal C(r_1)}|\omega|^\frac 52dx\Big)^\frac 9{40}\Big(\int\limits_{\mathcal C(r_1)}|\psi^2\omega_N|^\frac {s'115}{2(40-s'9)}dx\Big)^\frac {40-s'9}{40s'}.$$
Applying H\"older inequality one more time, we find 
$$\int\limits^0_{-r_1^2} \Big(\int\limits_{\mathcal C(r_1)} |\omega|^\frac {s'9}{16} |\psi^2\omega_N|^\frac {s'23}{16}dx\Big)^\frac 1{s'}\leq \Big(\int\limits_{Q(r_1)}|\omega|^\frac 52dz\Big)^\frac 9{40}\times$$$$\times\Big(\int\limits^0_{-r_1^2}\Big(\int\limits_{\mathcal C(r_1)}|\psi^2\omega_N|^\frac {s'115}{2(40-s'9)}dx\Big)^\frac {40-s'9}{31s'}dt\Big)^\frac {31}{40}.$$
Let us introduce numbers 
$$p=\frac {s'115}{2(40-s'9)}=\frac {23}4, \qquad q=p\frac {40-s'9}{31s'}=\frac {115}{62}.$$
It is easy to check that $3/p+2/q-3/2=1/10$ and thus 
 the known multiplicative inequality implies the bound

$$\|\psi^2\omega_N\|^\frac {23}{16}_{p,q,Q(r_1)}:=\Big(\int\limits^0_{-r_1^2}\Big(\int\limits_{\mathcal C(r_1)}|\psi^2\omega_N|^\frac {s'115}{2(40-s'9)}dx\Big)^\frac {40-s'9}{31s'}dt\Big)^\frac {31}{40}\leq 
$$
$$\leq c (r_1^\frac {1}{10}|\psi^2\omega_N|_{2,Q(r_1)} )^\frac {23}{16}.   $$
So, the final estimate for the $J_3$ has the form
$$J_3\leq \frac c{r_1-r}\Big(\int\limits_{\mathcal C(r_1)}|b|^sdx\Big)^\frac 1s\Big(\int\limits_{Q(r_1)}|\omega|^\frac 52dz\Big)^\frac 9{40} (r_1^\frac {1}{10}|\psi^2\omega_N|_{2,Q(r_1)} )^\frac {23}{16}\leq $$
$$\leq \frac c{r_1-r}r_1^\frac {23}{40}\Big(\int\limits_{Q(r_1)}|\omega|^\frac 52dz\Big)^\frac 9{40} (r_1^\frac {1}{10}|\psi^2\omega_N|_{2,Q(r_1)} )^\frac {23}{16}.$$
Now, our energy estimate can be re-written in the following way:
$$|\psi^2\omega_N|^2_{2,Q}\leq c\frac {r_1}{(r_1-r)^2}\Big(\int\limits_{Q(r_1)}|\omega|^\frac 52dz\Big)^\frac 45+$$
$$+\frac c{r_1-r}\Big(\int\limits_{Q_1(r_1)}|v|^\frac {10}3dz\Big)^\frac 3{10}\Big(\int\limits_{Q_1(r_1)}|\omega|^\frac {5}2dz\Big)^\frac 2{5}|\psi^2\omega_N|_{2,Q}+
$$
$$+\frac c{r_1-r}r_1^\frac {23}{40}\Big(\int\limits_{Q(r_1)}|\omega|^\frac 52dz\Big)^\frac 9{40} (r_1^\frac {1}{10}|\psi^2\omega_N|_{2,Q(r_1)} )^\frac {23}{16}.$$
After application of the Young inequality, we arrive at the estimate
$$|\psi^2\omega_N|^2_{2,Q}\leq \frac c{r_1-r}\Big(\Big(\frac 
{r_1}{r_1-r}\Big)^\frac {23}9+\Big(\int\limits_{(Q(r_1)}|v|^\frac {10}3dz\Big)^\frac 35\Big)\int\limits_{Q(r_1)}|\omega|^\frac 52dz\Big)^\frac 45,
$$
which, by the multiplicative inequality and by passing to the limit as $N\to\infty$, leads us to the final estimate of Step II:
\begin{equation}
	\label{basic}
	\Big(\int\limits_{Q(r)}|\omega|^\frac {10}3\Big)^\frac 3{10}\leq \frac c{\sqrt{r_1-r}}\Big(\Big(\frac 
{r_1}{r_1-r}\Big)^2+\Big(\int\limits_{(Q(r_1)}|v|^\frac {10}3dz\Big)^\frac3{10}\Big)\int\limits_{Q(r_1)}|\omega|^\frac 52dz\Big)^\frac 25 \leq
\end{equation}
$$\leq  \frac {cM}{\sqrt{r_1-r}}\Big(\frac 
{r_1}{r_1-r}\Big)^2\Big(\int\limits_{Q(r_1)}|\omega|^\frac 52dz\Big)^\frac 25, $$
where 
$$M=1+\Big(\int\limits_{(Q(r_1)}|v|^\frac {10}3dz\Big)^\frac3{10}.
$$

\textit{Step III}. For $k=2,3,...$,  set
$$m=m_k=\Big(\frac 43\Big)^k, \qquad r_1=r^{(k)}=\frac 12 +\frac 1{2^{k+1}},\, r=r^{(k+1)},\qquad Q_k=Q(r^{k}).
$$
Then, the following sequence of inequalities follows from the inequality (\ref{basic}):
$$\Big(\int\limits_{Q_{k+1}}|\sigma|^\frac {10m_k}3\Big)^\frac 3{10m_k}\leq $$\begin{equation}
	\label{sequence}\leq \Big(	\frac {cM}{\sqrt{r^{(k+1)}-r^{(k)}}}\Big(\frac 
{r^{(k+1)}}{r^{(k+1)}-r^{(k)}}\Big)^2\Big)^\frac 1{m_k}\Big(\int\limits_{Q_k}|\sigma|^\frac {5m_k}2dz\Big)^\frac 2{5m_k}.
\end{equation}
Observing that $10m_k/3=5m_{k+1}/2$ and applying Moser's arguments, we find
$$\sup\limits_{z\in Q(1/2)}|\sigma(z)| \leq c(M)\Big(\int\limits_{Q(5/8)}|\sigma|^\frac {5m_2}2dz\Big)^\frac 2{5m_2}.
$$
The right hand side of the above inequality 
is finite and  can be estimated by a constant depending only on $M$ as it has been shown in Step I.

So, finally, we conclude 
that
\begin{equation}
	\label{max1}
	\sup\limits_{z\in Q(1/2)}|\sigma(z)|= \sup\limits_{z\in Q(1/2)}|x'||v_\varphi(z)|\leq c(M).
\end{equation}

Now,  returning to the proof Theorem \ref{mainresult},  assume, in contrary,   that $z=0$ is a Type I blowup, i.e., 
\begin{equation}
\label{TypeIrestrict}
0<g=g(0)<\infty.	
\end{equation}
As it has been already noticed, (\ref{TypeIrestrict}) implies  that
\begin{equation}
\label{conseq}
G=G(0)<\infty.	
\end{equation}
So, without lose of generality, we may assume
that
\begin{equation}\label{scale-invatiant}
0<L_0=\sup\limits_{0<r<1}A(r)+\sup\limits_{0<r<1}C(r)+\sup\limits_{0<r<1}E(r)+	\sup\limits_{0<r<1}D(r)<\infty.
\end{equation}
Here, for example,
$$A(r)=A(0,r)=\sup\limits_{-r^2<t<0}\frac 1r\int\limits_{\mathcal C(r)}|v(x,t)|^2dx$$
and so on.



Now, we can  rescale our function $v$ and $q$ around the origin as it has been described, for example,  in \cite{SerShi2018}, see Theorem 3.5 there. 
Let $\lambda_k\to 0$ be a sequence and let
$$ u^k(y,s)=\lambda_kv(x,t),\quad p^k(y,s)=\lambda_k^2q(x,t),$$
where $x=\lambda_ky$ and $t=\lambda_k^2s$. 
Passing $k\to\infty$,  we can find limit functions $u$ and $p$ of sequences $u^k$ and $p^k$ that have the properties $(\mathcal A)$:

\noindent
(i) $u$ is a local energy ancient solution in $Q_-=\mathbb R^3\times]-\infty,0[$, i.e., the pair $u$ and $p$ is a suitable weak solution in $Q(R)$ for any $R>0$;

\noindent
(ii) $u$ and $p$ is axially symmetric solution to the Navier-Stokes equations in $Q_-$;

\noindent
(iii) for any $R>0$,
$$A(u,R)+E(u,R)+C(u,R)+D(p,R)\leq L_0<\infty.
$$

In addition, as it follows 
\
from (\ref{max1}), 
\begin{equation}\label{boundGamma}
	\Gamma=\varrho u_\varphi\in L_\infty(Q_-)
\end{equation}
and the velocity $u$ is not trivial in the sense 
\begin{equation}
\label{nontrivial}	
\frac 1{a^2}\int\limits_{Q(a)}|u|^3dz\geq \varepsilon(L_0)>0\end{equation}
for all $0<a\leq1$, see also \cite{Seregin2016}.

Since $u$ and $p$ is an axially symmetric suitable weak solution, there exists a closed set $S^\Gamma$ in $Q_-$, whose 1D-parabolic measure $\mathbb R^3\times \mathbb R$ is equal to zero and $x'=0$ for any $z=(x',x_3,t)\in S^\Gamma$, such that any spatial derivative of $u$ (and thus of $\Gamma$) is H\"older continuous in $Q_-\setminus S^\Gamma$. 
 In particular, for all $t<0$, the function $\Gamma(0,x_3,t)=0$ for all $x_3\in \mathbb R\setminus S^\Gamma_t$, where $S^\Gamma_t$ is a closed set of  measure zero in $\mathbb R$.
 
 Moreover, by the known multiplicative inequalities and by $(\mathcal A)_{iii}$, we have
 \begin{equation}
	\label{NazCond}	
\sup\limits_{R>0}
\frac 1{R^\frac 12}\Big(\int\limits_{-R^2}^0\Big(\int\limits_{\mathcal C(R)}|u|^3dx\Big)^\frac 43dt\Big)^\frac 14\leq c(L_0). \end{equation}
 And the function $\Gamma$ satisfies the following heat equation with the drift
\begin{equation}\label{GammaEquation}
\partial_t\Gamma+\Big(u+2\frac {(x',0)}{|x'|^2}\Big)\cdot\nabla \Gamma-\Delta \Gamma=0\end{equation}
in $\mathbb R^3\setminus \{x'=0\}\times ]-\infty,0[$.

It has been shown in \cite{NazUr} that if $\Gamma$ is a Lipschitz continuous function in $Q_-$ then, under the above conditions, it must be identically equal to zero. This could reduce our problem to the case with no drift, i.e.,  $u_\varphi=0$. 
Now, the  task is to show that the above Liouville type theorem remains to be true in our situation as well.

First, we need to understand  differentiability properties of the function $\Gamma$   in the domain $|x'|>0$.
To this end, we observe that
$$|\partial_t\Gamma(z)-\Delta \Gamma(z)|\leq (\sup\limits_{z=(x,t)\in P(\delta,R;R)\times ]-R^2,0[}|u(z)|+2/\delta)|\nabla \Gamma(z)|
$$ for any $0<\delta<R$, where $P(a,b;h)=\{x: \,a<|x'|<b,\,|x_3|<h\}$.
Since $u$ is axially symmetric, the first factor on the right hand side is finite. This, by iteration, yields
$$\Gamma\in W^{2,1}_p(P(\delta,R;R)\times ]-R^2,0[)$$
for any $0<\delta <R<\infty$ and for any finite $p\geq 2$. 

We also know that (it follows from the partial regularity theory) that, for any $t<0$,
\begin{equation}
	\label{limit}
	\Gamma(x',x_3,t)\to 0 \quad\mbox{as}\quad |x'|\to 0
\end{equation}
for all $x_3\in\mathbb R^3\setminus S^\Gamma_t$. %

Now, define the  class $\mathcal V$ of functions $\pi:Q_-\to \mathbb R$ possessing the  properties:

\noindent
(i) there exists a closed set $S^\pi$ in $Q_-$, whose 1D-parabolic measure $\mathbb R^3\times \mathbb R$ is equal to zero and $x'=0$ for any $z=(x',x_3,t)\in S^\pi$, such that any spatial derivatives is H\"older continuous in $Q_-\setminus S^\pi$; 

\noindent
(ii) $$\pi\in W^{2,1}_2(P(\delta,R;R)\times ]-R^2,0[)\cap L_\infty(Q_-)$$
for any $<\delta <R<\infty$.


Now, we can state an analog  of Lemma 4.2 of \cite{NazUr} for the class $\mathcal V$.
\begin{lemma}
	\label{analog}
	Let $u$ and $p$ have all the properties  $(\mathcal A)$. Let $\pi\in\mathcal V$ be a non-negative function satisfying the following conditions:
\begin{equation}\label{GammaInequality}
\partial_t\pi+\Big(u+2\frac {(x',0)}{|x'|^2}\Big)\cdot\nabla \pi-\Delta \pi\leq 0\end{equation}
in $\mathbb R^3\setminus \{x'=0\}\times ]-\infty,0[$,
\begin{equation}
	\label{1st}
	\pi(0,x_3,t)\geq k
	\end{equation}
	for all $t\in ]-R^2,0[$, $x_3\in S^\pi_t\cap ]-2R,2R[$ and for some $k>0$,	
	\begin{equation}
		\label{2nd}\pi\leq Mk\end{equation}
in $B(2R)\times ]-R^2,0[$ and for some $M\geq 1$.		Then
		\begin{equation}
			\label{consequence}\pi\geq \beta k
		\end{equation}
in $Q(R/2)$, where $\beta$ depends on $M$ and $N_R=R^{-1/2}\|u\|_{3,4}(Q(R))$ only.
	\end{lemma} 
\begin{proof} The proof of the lemma is based on the inequality:
$$\int\limits_{B(R)\times ]-R^2,-3R^2/4[}(\partial_t\pi\eta+\nabla \pi\cdot\nabla \eta-(u+2x'/|x'|^2)\cdot\nabla \eta\pi)dxdt\geq $$$$\geq4\pi_0
\int\limits_{B(R)\times ]-R^2,-3R^2/4[}\pi\eta dx_3dt
$$
for any  test function $\eta$ being equal to zero near spatial boundary. Here, $\pi_0=3.14...$. 
One can easily verify that the above inequality is still true for functions $\pi$ from the class $\mathcal V$.

Selecting test function in couple of special ways and repeating the arguments of the paper  \cite{NazUr}, 
we complete the proof. 
\end{proof}
Now, using standard arguments, we easily deduce from Lemma \ref{analog} that $\Gamma=0$ and, therefore, $u_\varphi=0$, i.e., $ u=u_\varrho e_\varrho+u_3e_3$.
Let us notice 
that 
any axially symmetric suitable weak solution with no swirl, i.e., $u_\varphi=0$, is smooth in the sense that any spatial derivative of $u$ is H\"older continuous in space-time (it can be done by considering a problem for $\eta=\omega_\varphi/\varrho$, where $\omega_\varphi=u_{3,\varrho}-u_{\varrho,3}$, and reduction of it to spatial dimension 5, see, for example,  \cite{KangK2004}).  In particular, the function $u$ is a continuous function  in $\overline Q(R)$ for any $R>0$. The latter contradicts restriction (\ref{nontrivial}) for sufficiently small  $a$.  Hence, the origin cannot be a Type I singularity. So,  as, by assumptions, $z=0$ is a singular point, it  should be Type II blowup.
\end{proof}



\setcounter{equation}{0}
\section{Axially symmetric solutions to the Navier-Stokes equations that are locally in critical spaces}

 Theorem \ref{mainresult} reads: axially symmetric solutions to the Navier-Stokes equations have no Type I blowups.   Therefore, any additional assumptions that exclude singularities of Type II, could be, in fact, sufficient conditions of regularity for axially symmetric solutions. For example, boundedness of norms in certain critical spaces is exactly such a case. It is known that solutions, belonging to critical space like $H^\frac 12$ or $L^3$, are smooth, see \cite{ESS2003}. As to weaker critical spaces like $L^{3,\infty}$, $BMO^{-1}$, or $\dot B^{-1}_{\infty,\infty}$, it is still unknown whether solutions, belonging to  these spaces, have no blowups. However, we can show  that the assumption of boundedness of suitable weak solutions in the above critical spacess excludes Type II singularities, see papers \cite{LeiZhang2011},  \cite{Seregin2011} and \cite{SerZhou2018}. Such an assumption and Theorem \ref{mainresult}  would imply  regularity of axially symmetric solutions.
 

Let us consider first the case of the space $L_\infty(-1,0;BMO^{-1}(\mathcal C))$. Here, in 3D case, $v={\rm rot}\,\omega$ with divergence free vector field $\omega\in L_\infty(-1,0;BMO(Q;\mathbb R^3))$ for which 
$$
\|v(\cdot,t)\|_{BMO^{-1}(\mathcal C)}:=     \|\omega(\cdot,t)\|_{BMO(\mathcal C)}=$$$$=\sup\Big\{\frac 1{|\mathcal C(r)|}\int\limits_{\mathcal C(x,r)}|\omega(y,t)-[\omega]_{\mathcal C(x,r)}(t)|dy:\,\,\mathcal C(x,r)\subset \mathcal C\Big\}
$$
$$\leq M<\infty$$
for a.a. $t\in ]-1,0[$. Then, according Proposition 1.1 of \cite{Seregin2011},
\begin{equation}\label{uniform1}
	\sup\limits_{z_0\in\overline Q(1/2),0<r<1/4} A(z_0,r)+D(z_0,r)+E(z_0,r)+C(z_0,r)<\infty.\end{equation}
and, hence, $z=0$ is a regular point, see Theorem \ref{mainresult}.

In the second case, $v\in L_\infty(0,T;\dot B^{-1}_{\infty,\infty}(\mathbb R^3))$. Here, there is a problem with localisation in the space $\dot B^{-1}_{\infty,\infty}(\mathbb R^3)$. To avoid unessential issues, we shall consider the problem of regularity in the whole space, i.e., in $Q_T=\mathbb R^3\times ]0,T[$. Let $v\in L_\infty(0,T;\dot B^{-1}_{\infty,\infty}(\mathbb R^3))$ be a suitable weak solution to the Navier-Stokes equations in $Q_T$. By Theorem 1.2 of \cite{SerZhou2018}, for $z_0\in \mathbb R^3\times ]0,T]$, we have the estimate
$$\sup\limits_{0<r<r_0}G(z_0,R)\leq c_0(r_0^\frac 12+\|v\|_{L_\infty(0,T;\dot B^{-1}_{\infty,\infty}(\mathbb R^3))}^2+\|v\|_{L_\infty(0,T;\dot B^{-1}_{\infty,\infty}(\mathbb R^3))}^6),$$
where $r_0\leq \frac 12\min\{1,t_0\}$ and $c_0$ depends $C(z_0,1)$ and $D(z_0,1)$. And again by Theorem \ref{mainresult}, such solution should be regular, for example, H\"older continuous.

The last case, $v\in L_\infty(-1,0;L^{3,\infty}(\mathcal C))$, is the easiest one as directly from H\"older inequality it follows that
$$\frac 1r\int\limits_{\mathcal C(x,r)}|v(y,t)|^2dy\leq c\|v(\cdot,t)\|^2_{L^{3,\infty}(\mathcal C(x,r))}
$$ and thus 
$$\sup\limits_{z_0\in \overline Q(0,1/4)}\sup\limits_{0<r<1/4} A(z_0,r)\leq c\|v(\cdot,t)\|^2_{L^{3,\infty}(\mathcal C)}.$$
Then, an estimate of type (\ref{uniform1}) can be deduced from \cite{Seregin2006}, see Lemma 1.8, with the same conclusion, by Theorem \ref{mainresult},  that $z=0$ is a regular point.

{\bf Acknowledgement} The work is supported by the grant RFBR 20-01-00397.

\end{document}